\documentclass{amsart}
\usepackage{amsmath}
\usepackage{amssymb}
\usepackage{amsthm}
\usepackage{amscd}
\usepackage{amsfonts}
\usepackage{graphicx}%
\usepackage{fancyhdr}
\usepackage[all]{xy}
\usepackage[ mathscr]{euscript}
\usepackage{amssymb}
\usepackage{amsmath}
\usepackage{amsfonts}
\usepackage{graphicx}
\usepackage{natbib}
\bibpunct{[}{]}{;}{a}{,}{,}
\theoremstyle{plain} \numberwithin{equation}{section}
\newtheorem{theorem}{Theorem}[section]
\newtheorem{corollary}[theorem]{Corollary}
\newtheorem{conjecture}[theorem]{Conjecture}
\newtheorem{lemma}[theorem]{Lemma}
\newtheorem{proposition}[theorem]{Proposition}
\theoremstyle{definition}
\newtheorem{definition}[theorem]{Definition}

\newcommand{\Z}{\mathbb{Z}}
\title{Finite Groups With Maximal Normalizers I}
\author{Joseph P. Bohanon}
\thanks{The author would like to thank his advisor, John Shareshian, as well as Yakov Berkovich and Eamonn O'Brien for all of their helpful comments on this paper.}

\begin{document}
\maketitle
\begin{abstract}
We examine $p$-groups with the property that every non-normal subgroup has a normalizer which is a maximal subgroup.  In particular we show that for such a $p$-group $G$, when $p=2$, the center of $G$ has index at most 16 and when $p$ is odd the center of $G$ has index at most $p^3$.
\end{abstract}

\section{Introduction}
In this paper, and the one to follow, we will examine groups in which the normalizer of every non-normal subgroup is a maximal subgroup.  We call such a group an MN-group.  This paper will examine $p$-groups.  The second paper will examine all other MN-groups.  The principal result of this paper is the following:
\begin{theorem} \label{mainresult} Let $G$ be a $p$-group in which every non-normal subgroup has $p$ conjugates. Then \begin{equation*}
[G:Z(G)] \le \begin{cases} 16 \text{ if } p=2, \\ p^3 \text{ if } p>2. \end{cases} \end{equation*}
\end{theorem}
If $H \le G$ has at most $p$ conjugates, then $[G:N_G(H)] \le p$, therefore $G$ is an MN-group.  Also, as $\Phi(G)=G' G^p$, every $p$-th power must normalize every subgroup of $G$.
\\ \indent In Section 2 we provide some of the background material and preliminary results.  In Section 3, we examine $p$-groups with element breadth 1.  In Section 4, we examine $2$-groups with element breadth 2.  In Section 5, we examine $p$-groups with element breadth 2 for odd $p$.  In Section 6, we make a conjecture that would generalize \ref{mainresult} to groups with subgroup breadth $k$.  In Section 7, we provide details about our use of GAP to solve this problem.
\section{Preliminary Results and Definitions}
All groups in this paper will be finite.  Recall that a group is \textit{Hamiltonian} if every subgroup is normal.  The following theorem of Dedekind is the starting point for this paper:
\begin{theorem} Let $G$ be a Hamiltonian group.  Then $G$ is either abelian, or $G \cong Q_8 \times (\Z_2)^n \times A$ where $A$ is an abelian group of odd order.
\end{theorem}
\begin{proof}
See \cite{robinson}.
\end{proof}
We begin with some definitions:
\begin{definition} The \textit{element breadth} of an element $x$ of a $p$-group $G$, $\mathrm{ebr}(x)$, is defined to be the integer such that $[G:C_G(x)]=p^{\mathrm{ebr}(x)}$.  The element breadth of
$G$, $\mathrm{ebr}(G)$, is the maximum value that $\mathrm{ebr}(x)$ takes over all the elements
of $G$. \end{definition}
It should be noted that ``element breadth" is a non-standard term (\textit{breadth} is the standard term).  We use the term to distinguish between the previous definition and the following one.
\begin{definition} The \textit{subgroup breadth} of a subgroup $H$ in a $p$-group $G$, $\mathrm{sbr}(H)$, is defined to be the integer such that $[G:N_G(H)]=p^{\mathrm{sbr}(H)}$.  The subgroup breadth of $G$, $\mathrm{sbr}(G)$, is the maximum value that $\mathrm{sbr}(H)$ takes over all the subgroups of $G$.  The \textit{cyclic breadth} of $G$, $\mathrm{cbr}(G)$ is the maximum value that $\mathrm{sbr}(H)$ takes over all cyclic subgroups of $G$.
\end{definition}
As an example, consider an extra-special 2-group $G$ of plus-type of order $2^{2m+1}$.  Clearly every non-normal subgroup must be elementary abelian of order at most $2^m$.  Let $H$ be a non-normal subgroup of order $2^m$ generated by non-central involutions $x_1,\cdots,x_m$.  Then $C_G(H)=\displaystyle \bigcap_{i=1}^m C_G(x_i)$ has index at most $2^m$ so that $\mathrm{sbr}(G) \le m$.  As $G$ is of plus-type, we have generators $i_1,j_1,\cdots,i_m,j_m$ such that $\langle i_k,j_k \rangle \cong Q_8$ for $1 \le k \le m$, each such $Q_8$ commutes with every other $Q_8$ and the $-1$ from each of the $Q_8$'s have been identified.  In particular, a non-central involution of $G$ is a product of involutions of the form $a_s b_t$ where $s \ne t$, $a_s$ is one of $i_s,j_s$ or $k_s$ and $b_s$ is one of $i_t,j_t$ or $k_t$.  Consider the subgroup $K=\langle i_1 i_2,j_1 j_2 \cdots,i_{2m-1} i_{2m},j_{2m-1} j_{2m} \rangle$.  This group is elementary abelian and is clearly non-normal.  Moreover, if some element $g$ normalizes $K$ it cannot conjugate the element $i_s i_t$ to $j_s j_t$ or $k_s k_t$ therefore $N_G(K)=C_G(K)$.  It is easily seen that the centralizers of the generators of $K$ are mutually distinct, therefore $[G:N_G(K)]=p^m$ so that $\mathrm{sbr}(G) \ge m$, therefore $\mathrm{sbr}(G)=m$.  It can similarly be shown that an extra-special group of minus-type has subgroup breadth $m-1$.  Note also that every cyclic subgroup of order more than 4 must contain the center, therefore must be normal, therefore $\mathrm{cbr}(G)=\mathrm{ebr}(G)=1$.  Therefore cyclic breadth has no influence on subgroup breadth.
\\ \indent We will need use of GAP \cite{GAP4} as well. In particular, we will use the notation $[n,m]$ to be group $m$ of order $n$ in the Small Group Library \cite{smallgroups}.  (It should be noted that the Small Group Library is not specific to GAP.)  We access the library and use the GAP functions \textit{ConjugacyClasses} and \textit{ConjugacyClassesSubgroups} to determine the element and subgroup breadths of a $p$-group.  We will repeatedly use the results of section 7 to show that certain groups involved in minimal counterexamples do not have subgroup breadth 1. Finally, it should be noted that some of the groups that we must construct in GAP are constructed as finitely-presented groups, therefore we use the GAP function \textit{IsomorphismPcGroup} to speed up necessary computations over what can be expected in a finitely presented group.
\\ \indent For $p=2$, \ref{mainresult} was first posed as a conjecture in \cite{knoche}.  In that paper it was shown that condition (TC), that is $\mathrm{sbr}(G)=1$, is equivalent to each of the following two conditions:
\\\\ (CO) The core $H_G$ of every subgroup $H$ of $G$ ``requires" all the conjugates of $H$, in the sense that the intersection of a proper subset of the set of distinct conjugates of $H$ properly contains $H_G$.
\\\\ (NC) The normal closure $H^G$ of every subgroup $H$ of $G$ ``requires" all the conjugates of $H$, in the sense that the subgroup generated by a proper subset of the set of distinct conjugates of $H$ is a proper subgroup of $H^G$.
\\\\ \indent The question of whether such a bound on the index of the center exists or not is found in \cite{bj2} as suggested research problem 830.  In \cite{mac} it is shown that \begin{equation*} [G:Z(G)] < p^{81 \mathrm{sbr}(G) (\log_2 (p^{\mathrm{sbr}(G)^2}))}, \end{equation*} therefore the results in this paper are a great improvement over this result in the case that $\mathrm{sbr}(G)=1$.
\\ \indent The central product of $Q_8$ and $D_8$ has subgroup breadth 1 and a center of
index 16, so the bound is the best possible when $p=2$.  For $p=3$, the group 
$$ \langle a,b,c \mid a^{p^3}=b^p=c^p=1,[a,b]=a^{p^2},[a,c]=a^{p^2} b,[b,c]=1 \rangle$$
shows the bound is sharp.  For $p>3$ the group $$ \langle a,b,c \mid a^{p^2}=b^p=c^p=1,[a,b]=a^{p},[a,c]=a^{p} b,[b,c]=1 \rangle$$ shows the bound is sharp.
\\ \indent The following result bounding the element breadth of a $p$-group with given cyclic breadth \cite{csw}.
\begin{proposition} If $G$ is a $p$-group, then
\begin{equation*} \mathrm{ebr}(G) \le \begin{cases} 2 \mathrm{cbr}(G)+1 & \text{if } p=2,
\\ 2 \mathrm{cbr}(G) & \text{if } p>2. \end{cases}
\end{equation*}
\end{proposition}

In particular, this says that $\mathrm{ebr}(G) \le 3$ when $\mathrm{sbr}(G)=1$.

\begin{proposition} \cite{knoche} \label{br1} If $G$ is a $p$-group, then $\mathrm{ebr}(G)=1$ if and only if $|G'|=p$.
\end{proposition}

\begin{proposition} \cite{ps} \label{br2} If $G$ is a $p$-group, then $\mathrm{ebr}(G)=2$ if and only if one of the following
holds:
\begin{enumerate}
\item $|G'|=p^2$ or
\item $[G:Z(G)]=p^3$ and $|G'|=p^3$.
\end{enumerate}
\end{proposition}

The following can be found in \cite{blackburnint}:
\begin{proposition} \label{blackburnint} Let $G$ be a $p$-group in which the intersection of all the non-normal subgroups is non-trivial.  Then $p=2$ and $G$ is isomorphic to one of the following:
\begin{enumerate}
\item $C \cong Q_8 \times \Z_4 \times (\Z_2)^k$.
\item $C \cong Q_8 \times Q_8 \times (\Z_2)^k$.
Here $Q_8 \times Q_8$ does not have subgroup breadth 1.
\item $C=\langle g,A \rangle$ with $A$ abelian but not
elementary abelian, $1 \ne g^2 \in A$ and $a^g=a^{-1}$.
\end{enumerate}
\end{proposition}

The first reduction we make is the following.

\begin{theorem} If $\mathrm{sbr}(G)=1$ then $\mathrm{ebr}(G) \le 2$. \end{theorem}

We give several results about metacyclic groups and determine
which metacyclic groups have subgroup breadth 1.  The following may be found in \cite{king}.

\begin{lemma} \label{mcpres} If $G$ is a metacyclic $p$-group, then $G$ has a presentation of the form
$$\langle a,b \mid a^{p^m}=1,b^{p^n}=a^k, a^b=a^r \rangle$$ where $$m,n \ge 0, 0<r, k<p^m, p^m \mid k (r-1), \text{ and } p^m \mid r^{p^n}-1.$$  Up to
picking different generators, we may assume that $k=0$ or $p^j$.
\end{lemma}
\begin{proof}
The first statement is proved in \cite{king}.  For the second statement, if $k \ne 0$, we may obtain the result by replacing $a$ or $b$ with powers of $a$ or $b$.
\end{proof}
\begin{corollary} \label{metacor} Let $G$ be a 2-group as in \ref{mcpres} with $\mathrm{ebr}(G)$=2 and $\mathrm{sbr}(G)$=1. Then $|\langle a \rangle \cap \langle b \rangle| \le 4$.
\end{corollary}
\begin{proof}
By \ref{br2} we know that $G'$ has size either 4 or 8 and we can write $$G=\langle a,b \mid a^{2^m}=1,b^{2^n}=a^{2^{m-\ell}},[a,b]=a^{m-2} \rangle $$ or $$G=\langle a,b \mid a^{2^m}=1,b^{2^n}=a^{2^{m-\ell}},[a,b]=a^{m-3} \rangle .$$
We assume that $\ell \ge 3$.
Consider the first case.  Note that we have that $\langle b \rangle$ is also normal in $G$, so without loss of generality, we may assume that $|a| \ge |b|$.  Note that $|G|=2^{m+n}$.  Now, $a^4$ and $b^4$ are central, therefore if $m \ge 5$, $\langle b^{-4} a^{2^{m-3}} \rangle \lhd G$.  Quotienting out by this subgroup, we get a corresponding metacyclic group with $n=2$ and $\ell=3$.  When $m=5$, consider the subgroup $\langle a^3 b \rangle$ and when $m>5$ consider the subgroup $\langle a^{-2^{n-4}} b \rangle$.  These groups both have order 4, and intersect $\langle a \rangle$ trivially, therefore by picking different generators, we may assume that $G$ is a split extension $\Z_{2^m} \rtimes \Z_4$.  Letting $c$ be the generator of $\Z_4$, we still have that $[a,c]=a^{\pm 2^{m-2}}$.  In this group $\langle c \rangle$ has four conjugates.  Therefore, we may assume that $m=4$ (when $m=3$, $G$ is clearly a split extension), hence $\ell=3$.  It is easily verified, in GAP, that the groups with the given presentations for $n=2,3$ or 4, have element breadth 1.
\\ \indent Consider the second case.  Again, $\langle b \rangle$ is normal, so we assume $|a| \ge |b|$.  Also, $\langle b^{-8} a^{2^{m-3}} \rangle$ is normal in $G$ if $m \ge 6$.  When $m=6$ consider the subgroup $\langle a^3 b \rangle$ and when $m > 6$ consider the subgroup $\langle a^{-2^{m-4}} b \rangle$.  As before, these groups intersect $\langle a \rangle$ trivially, therefore, we can write $G$ as a semi-direct product $\Z_{2^m} \rtimes \Z_8$.  However, if $c$ is a generator of $\Z_8$, then $\langle c \rangle$ has more than two conjugates.  Therefore $m<6$.  If $m=5$, suppose $\ell=3$.  When $n=1$, $G$ has element breadth 1.  When $n=2$ we get [64,28] and when $n=3$ we get [128,130] neither of which have subgroup breadth 1 by \ref{64grps} and \ref{128grps}.  It is easily checked that when $\ell=4$ that $G$ has element breadth 1 in all possible cases.  When $m=4$, similarly, all remaining groups have element breadth 1.  When $m=3$, $G$ is cyclic, therefore we are done.
\end{proof}
Now, we prove some structure theorems regarding metacyclic 2-groups with
subgroup breadth 1.
\begin{lemma} Let $G$ be a metacyclic $p$-group with the presentation as
in \ref{mcpres}.  Then we have $$(ba^i)^j=b^j a^{i
(1+r+r^2+\cdots+r^{j-1})}$$
\end{lemma}
\begin{proof} By induction.  When $j=1$ this is obvious.  Now
suppose we have the statement for $j-1$.  Then
$$(ba^i)^j=(ba^i)^{j-1} ba^i=b^{j-1} a^{i (1+r+r^2+\cdots+r^{j-2})} b
a^i=$$ $$b^j a^{r i (1+r+r^2+\cdots+r^{j-2})} a^i=b^j a^{i
(1+r+r^2+\cdots+r^{j-1})}.$$ \end{proof}
We prove some results about 2-groups.
\begin{theorem} If a 2-group $G$ is metacyclic, non-abelian and
$\mathrm{sbr}(G)=\mathrm{ebr}(G)=1$ then with the notation from \ref{mcpres} we have
$r=2^{m-1}+1$.
\end{theorem}
\begin{proof} We have $[a,b]=a^{r-1}$.  Since $\mathrm{ebr}(G)=1$, by \ref{br1} we must
have that $a^{r-1}$ has order 2.  This says that $2r \equiv 2 \pmod{2^m}$, which says that $r=2^{m-1}+1$ $(r \ne 1$ since $G$ is
non-abelian).
\end{proof}
\begin{theorem} \label{br2meta} If a 2-group $G$ is metacyclic,
$\mathrm{sbr}(G)=1$, and $\mathrm{ebr}(G)=2$ then with the notation as above, there is some $n$ such that
$$G=\langle a,b \mid a^{8}=b^{2^n}=1,a^4=b^{2^{n-1}},[a,b]=a^{\pm 2}
\rangle.$$
\end{theorem}
\begin{proof}
Note first that $G'=\langle [a,b] \rangle$.  This is true because clearly $\langle a^2 \rangle \lhd G$ and $G/\langle a^2 \rangle$ is abelian, therefore $G' \le \langle a^2 \rangle$ and no proper subgroup can have an abelian quotient.  Suppose that $\langle a \rangle
\cap \langle b \rangle=1$.  If $a$ has order 2, then $a$ and $b$
commute and clearly $G$ is abelian. Consider the subgroup $\langle b
\rangle$.  We have $\langle b
\rangle^a=\langle ba^{-r+1} \rangle$ and $\langle b
\rangle^{a^2}=\langle ba^{-2r+2} \rangle$.  Hence we must have
that $\langle b \rangle=\langle
ba^{-2r+2} \rangle$.  In this case, $r=2^{m-1}+1$ so $G$
has element breadth 1 by \ref{br1}, because $G'$ has order 2.
This shows that $\langle a \rangle$ and
$\langle b \rangle$ intersect non-trivially.  We showed in \ref{metacor} that we cannot have that $|G'|=8$.  Therefore $G'$ must have order 4,
and we may assume that $r=2^{m-2}+1$ or $3 \cdot 2^{m-2}+1$.  Suppose first that $\langle b \rangle=\langle b^a \rangle=\langle b
a^{-r+1} \rangle$.  Then $a^{2^{m-2}}=b^{\pm 2^{n}}$ by \ref{metacor}.  We must
examine these groups for some small values of $m$ and $n$ first.
Suppose that $m=4$ so that $r=5$ or 13 and we have $a^4=b^{\pm
2^n}$.  When $n=2$, $G$=[64,28] and when $n=3, G=$[128,130], neither
of which have subgroup breadth 1 by \ref{64grps} and \ref{128grps}.  Both of these groups also have a unique normal subgroup
of order $2$ generated by $a^8$, hence no quotient groups satisfy our hypothesis that $a$ has order more than 8.
 So we may assume that $n \ge 4$.
Now, since $H=H^{b^2}$ we have that $a^8=b^{2^{n-1}}$ is a power of
$b^{2^{n-2}} a^{\mp 1}$ which has order 4 (throughout this proof the use of $\pm$ and $\mp$ indicates that one must consistently pick the sign on the top or the sign on the bottom). Since $a^8$ has order 2,
we must have $(b^{2^{n-2}} a^{\mp 1})^2=b^{2^{n-1}} a^{\mp 2}=a^8$.
This says that $a^2$ is a power of $b$.  However $(a^2)^b=a^{10}$
which says that $a^2=a^{10}$ and $a$ has order 8, a contradiction.
Therefore $m \ne 4$.
\\ \indent Now, assume that $m>4$. We proceed by
induction on the size of $G$. Given the relations for $G$ it is easily verified that $a^4$ and $b^4$ are central.  Moreover, $\langle a^{2^{m-2}} \rangle$ and $\langle a^{2^{m-3}} b^{2^{n-1}} \rangle$ are two central subgroups of order 4, therefore $Z(G)$ is not cyclic.  Therefore, there is some involution $z$ which is not a commutator.  We get that $\overline{G}=G/\langle z \rangle$ is a metacyclic group with element breadth 2.
By induction the order of $\overline{a}$ in this
quotient group must be 8.  This says that $a$ has order 8 or 16,
contradicting our assumption.  Hence we must have that $n=1$ or 2.
If $n=1$, then $b^2$ commutes with $a$.  However
$a^{b^2}=a^{2^{m-1}+1}$.  This says that $a$ has order $2^{m-1}$ and
we get the result by induction.  Now, we assume that $n=2$.  We
claim that this group is actually a split extension isomorphic to
$\Z_{2^m} \rtimes \Z_4$.  Let $\langle a \rangle$ be the subgroup
isomorphic to $\Z_{2^m}$.  Consider the subgroup $\langle
a^{2^{m-4}} b^{\mp 1} \rangle$.  When $m=5$ we again have that $G$ is a quotient of
[128,130],
so we assume that $m \ge 6$.  Then $a^{2^{m-4}}$ and $b$ commute,
hence $(a^{2^{m-4}} b^{\mp 1})^4=a^{2^{m-2}} b^{\mp 4}=1$.  This
contradicts our earlier statement that the generators of a metacyclic
group with element breadth 2 and subgroup breadth 1 generate
subgroups with non-trivial intersection.
\\ \indent Therefore, we may assume that $m=3$.  So $G$ has a presentation of
the form $$\langle a,b \mid a^8=1,b^{2^n}=a^{2^k},[a,b]=a^{\pm 2}
\rangle.$$  Note that such a group has a unique central involution, namely $a^4$, hence
no quotient groups still satisfy the requirement that $a$ has order 8.  If $k=3$, the subgroup generated by $b$ has four
conjugates.  If $k=1$, note that $1=[b^{2^n},b]=[a^2,b]=a^4$.  So we
may assume that $k=2$, which gives the lemma.
\end{proof}
By a result in \cite{csw}, if $G$ is a metacyclic 2-group, then $\mathrm{ebr}(G) \le \mathrm{sbr}(G)+1$, therefore there are no metacyclic groups with element breadth 3 and subgroup breadth 1.
\\ \indent We remark that in the element breadth 2 case, the conjugacy class of $b$ is the only
one with four elements.  Also the subgroup $\langle b \rangle$ has two conjugates.
\\ \indent Recall that the only non-trivial automorphisms of order 2 of a cyclic group of order $2^n$ send a generator $a$
to $a^{-1}$, $a^{2^{n-1}-1}$ or $a^{2^{n-1}+1}$.  We call these automorphisms \textit{dihedral}, \textit{semi-dihedral} and \textit{modular}, respectively.
We can use the above results to show the following.
\begin{theorem} If $G$ is a 2-group with $\mathrm{ebr}(G)=3$ then $\mathrm{sbr}(G) \ne 1$. \end{theorem}
\begin{proof}
Let $a$ be an element with 8 conjugates.  We aim to
show that $a$ has order 8.
Let $A=\langle a \rangle$ and $H=N_G(\langle a \rangle)$. Suppose that $a$ has order $2^n$
where $n>3$. Then Aut$(A) \cong \Z_{2^{n-2}} \times \Z_2$ where the first factor is the automorphism $a \to a^5$ and the second factor is the automorphism $a \to a^{-1}$ and we get a homomorphism $\phi:H \to \mathrm{Aut}(A)$.  If $\phi(H)$ contains a dihedral automorphism,
let $t$ be a pre-image.  Then consider $\langle a,t \rangle$ and the subgroup $K=\langle t \rangle$.  Then $K=K^{a^2}=\langle t a^{-4} \rangle$.  Therefore $a^4 \in \langle t \rangle$, so $a$ has order 16 by \ref{metacor}.  However, in this case $[a,t]$ has order 8, a contradiction.  Similarly, if $\phi(H)$ contains a semi-dihedral automorphism, we again get that $a^4 \in \langle t \rangle$.  Therefore $\phi(H)$ is cyclic.  Let $t$ be a
pre-image of a generator.  If $H=G$ then we get that $\langle a,t \rangle$ must contain all eight conjugates of $a$, however, this group is metacyclic, a contradiction.  Therefore $[G:H]=2$.  Let $C=C_G(a)$.  Then $|H/C|=4$.  If $H/C$ is cyclic, then let $t$ be the pre-image of a generator.  Again $\langle a,t \rangle$ is metacyclic and since it has element breadth 2, $a$ would have order 8.  However, this is impossible since $\mathrm{Aut}(A)$ has no elements of order 4.  Therefore $N/C \cong \Z_2 \times \Z_2$.  Let $b$ and $c$ be pre-images of the generators of this group.  We may assume that $b$ induces a dihedral automorphism and $c$ induces a semi-dihedral automorphism.  Now, $b^2$ and $c^4$ are central elements.  Then
$\langle a,b \rangle$ has element breadth at least 2, hence $a$ must
have order 8 and $a^4=b^{2^n}$ and similarly $a^4=c^{2^\ell}$. Suppose
first that $b$ does not have eight conjugates.  In $\langle a,b
\rangle$, $b$ already has four conjugates, hence $b^c$ must be one
of these, that is, $[b,c]=a^{2i}$ where $i=0,1,2$ or 3.  Note that
this implies that $b^2$ and $c^4$ are central elements.  Suppose
first that $m \ge 3$. We look at the subgroup $H=\langle
ab^{2^{m-2}} \rangle$.  Conjugating by $b$ and $c$ we get
$H^b=\langle a^{-1} b^{2^{m-2}} \rangle$ and $H^c=\langle a^3
b^{2^{m-2}} \rangle$ respectively.  If $H=H^b$ or $H=H^c$ then $a^2$ is a
power of $ab^{2^{m-2}}$ which has order 4.  This implies that $a$ is
a power of $b$, a contradiction.  If $H^b=H^c$ then $a^4$ is a power
of $ab^{2^{m-2}}$.  Since $a^4$ has order 2, this implies that
$a^4=a^2 b^{2^{m-1}}$ which is also impossible by the structure of a
metacyclic group with element breadth 2 and subgroup breadth 1. This
shows that if $b$ has 4 conjugates, it must have order at most 8.
Similarly, $c$ has order at most 16.  Note that since $\langle b \rangle$ has two conjugates in $\langle a,b \rangle$, we must have that $\langle a,b \rangle \lhd \langle a,b,c \rangle$.  Therefore $|\langle a,b,c \rangle| \le 256$ and has element breadth 3.  By \ref{256grps} there are no such groups.  This completes the proof.
\end{proof}
Next we prove two unpublished results of John Shareshian.
\begin{proposition} \label{johnd8} If $G$ is a 2-group with subgroup breadth 1 and
two involutions of $G$ do not commute, then the center of $G$ has
index at most 16.
\end{proposition}
\begin{proof} Let $G$ be a 2-group and let $s$ and $t$ be involutions such that $[s,t] \ne 1$.
 Then $\langle s,t \rangle$ is a dihedral group.
 Since $D_{2^n}$ does not have subgroup breadth 1 for $n \ge 4$,
 we clearly get $\langle s,t \rangle=D \cong D_8$.  Since $s$ and $t$
 already have two conjugates in $D$, we must have $D \lhd G$.
 Let $C=C_G(D)$.  Then $G=CD$ (see, for instance, 4.17 in \cite{suzuki2}).  Let $z$ generate $Z(D)$.
 We must have $C \cap D \le \langle z \rangle$.
 If $H<C$ is not normal in $C$, then $H$ must contain $z$
 (otherwise if $H \ne K$ and $H$ is conjugate to $K$, the four groups
 $\langle t,H \rangle, \langle t,K \rangle, \langle tz,H \rangle$ and
 $\langle tz,K \rangle$ are distinct and conjugate).
 If $C$ has no non-normal subgroup, then $C$ is Hamiltonian and it is
 straightforward that the proposition holds.  Otherwise, by \ref{blackburnint}, one of the following occurs:
\begin{enumerate}
\item $C \cong Q_8 \times \Z_4 \times (\Z_2)^k$.
It is straightforward to verify that the proposition holds in this case.
\item $C \cong Q_8 \times Q_8 \times (\Z_2)^k$.
Here $Q_8 \times Q_8$ does not have subgroup breadth 1.
\item $C=\langle g,A \rangle$ with $A$ abelian but not
elementary abelian, $1 \ne g^2 \in A$ and $a^g=a^{-1}$.
Since $g$ both centralizes and inverts $g^2$, $g$ has
order 4 and $g^2=z$.  Suppose some $a \in A$ has order
eight.  If $g^2 \ne a^4$ then the involutions $stg,stga^2,tsg$
and $tsga^2$ are distinct and conjugate, a contradiction.
Next suppose there is some $b \in A$ such that $b$ has order
4 and $b^2 \ne g^2$.  Then $g^b=gb^2$, so $stg,stgb^2,tsg$
and $tsgb^2$ are distinct and conjugate, a contradiction.
Therefore $A$ is the direct product of a cyclic group of
order 4 and an elementary abelian group.  Therefore $[C:Z(C)]=4$
and it is straightforward to verify the proposition.
\end{enumerate}  \end{proof}
\begin{theorem} \label{3inv} If $\mathrm{sbr}(G)$=1, $Z(G)$ is cyclic and all involutions of $G$ commute with each other, then $G$ contains at most three involutions.
\end{theorem}
\begin{proof}
Let $z$ be the unique element of order 2 in $Z(G)$ and let $t \ne z$ have order 2 in $G$.  Let $t' \ne t$ be a conjugate of $t$.  Since $[G:C_G(t)]=2$ we have $C_G(t) \lhd G$, so $C_G(t')=C_G(t)$.  Therefore if $g \in C_G(t)$, $(tt')^g=tt'$, while if $g \in G-C_G(t)$, then $(tt')^g=t't=tt'$.  Thus $tt' \in Z(G)$ and since $t$ and $t'$ commute, we have $|tt'|=2$, so $t'=tz$.  Suppose there is some involution $s$ in $G$ besides $t,t'$ and $z$.  Since $s$ and $t$ commute, $|st|=2$ and $st \notin \{z,t,t'\}$.  If $C_G(s)=C_G(t)$, let $g \in G-C_G(t)$.  Then $(st)^g=sztz=st$, so $C_G(st) \ne C_G(t)$.  Then $\langle s,t \rangle,\langle s,tz \rangle,\langle sz,t \rangle$ and $\langle sz,tz \rangle$ are distinct and conjugate, a contradiction.
\end{proof}
The result of Blackburn (\ref{blackburnint}) will be used multiple times.  We examine the groups in case 3 in more detail.  Since $g^2 \in A$ any element
of this group not in $A$ has the form $ga$ where $a \in A$. As mentioned above, $g$ must have order 4. Also, $g^2$ is
central.  Now, every involution in $A$ must centralize $g$, therefore so does $\Omega_1(A)$.  We must have that $[A:\Omega_1(A)] \le 4$.  If $A=\Omega_1(A)$ then $A$ is elementary abelian and $C$ is abelian.  Therefore, $A$ is either isomorphic to $\Z_8 \times (\Z_2)^n$ or $\Z_4 \times \Z_4 \times (\Z_2)^n$.  Consider the former case.  If $g^2$ does not lie in $\Z_8$, then letting $t$ be a generator of $\Z_8$, $\langle g,t \rangle$ is [32,14] which does not have subgroup breadth 1 by \ref{32grps}.  Otherwise, $\langle g,t \rangle$ is isomorphic to $Q_{16}$ and therefore contains a subgroup isomorphic to $Q_8$.  Consider the latter case.  If $g^2$ is not the square of an element of $A$ of order 4, then $C$ contains a subgroup isomorphic to $\Z_4 \times \Z_4$.  This group contains quotients isomorphic to $Q_8$ and $D_8$.  If $g^2$ is the square of an element of order 4 then $C$ contains a subgroup isomorphic to $Q_8$.
\\ \indent Suppose that we have some $a \in C$ of order 8. If $g^2 \ne
a^4$ then in $\langle g,a \rangle$ we have $\langle g \rangle,
\langle g a^2 \rangle, \langle g a^4 \rangle$ and $\langle g a^6
\rangle$ are conjugate. If $C \ge \Z_8 \times \Z_8$,
one of the generators of the direct factors of $\Z_8 \times \Z_8$ cannot be equal to $g^2$.
So we can
assume that $C \cong \Z_{2^i} \times (\Z_4)^m \times (\Z_2)^n$ or $C
\cong (\Z_4)^m \times (\Z_2)^n$.  Suppose now that $g^2$ is the
generator of a $\Z_2$. (Clearly if $g^2$ is not a square in $C$, we can always pick a presentation for $C$ such
that $g^2$ is such a generator.) Let $a$ and $b$ be elements such that
$|a|,|b|>2$ and $\langle a \rangle \cap \langle b \rangle=1$.
 Then in $\langle g,a,b \rangle$, consider the subgroup $H=\langle g
 \rangle$.  We have $H^a=\langle g a^2 \rangle$ and $H^b=\langle g
 b^2 \rangle$.  Since $g^2$ cannot be an element of one of the
 non-$\Z_2$ factors, these two groups are clearly distinct from $H$
 and from each other.  This says that we may assume that $C \cong \Z_{2^i}
 \times (\Z_2)^n$.  Also, if $g^2$ is the square of some element of
 order 4, then $B$ contains a quaternion subgroup.  However, only
 $Q_8$ and $Q_{16}$ have subgroup breadth 1, which says that $i=2$
 or $i=3$.  Hence, in this case, $B$ contains a subgroup isomorphic
 to $Q_8$.  Otherwise, since $g$ inverts the element of order $2^i$
 we have a subgroup isomorphic to $\langle a,b \mid
 a^{2^i}=b^4=1,[a,b]=a^2 \rangle$.  If $i>2$, the subgroup generated
 by $b$ has at least four conjugates, hence we may assume that $B
 \cong \Z_4 \rtimes \Z_4$ (where the action is irreducible).  By
 quotienting out by $b^2$ we get a quotient group isomorphic to $D_8$
 and by quotienting out by $a^2 b^2$ we get a quotient group isomorphic
 to $Q_8$.
\section{$p$-groups with element breadth 1}
By \ref{br1} we know that $|G'|=p$.  We first state two results from \cite{bj1}.  Recall that a minimal non-abelian group is a group all of whose proper subgroups are abelian.
\begin{theorem} \cite{bj1} (1.18a)
A minimal non-abelian $p$-group is isomorphic to one of the following:
\begin{enumerate}
\item $Q_8$,
\item $P_{i,j}=\langle a,b \mid a^{p^i}=b^{p^j}=1,[a,b]=a^{p^{i-1}} \rangle,i \ge 2,j \ge 1$,
\item $P_{i,1,k}=\langle a,b,c \mid a^{p^i}=b^p=c^{p^k},[a,b]=[b,c]=1,[a,c]=b \rangle,i+j>2.$
\end{enumerate}
\end{theorem}
\begin{theorem} \cite{bj1} (4.2) Let $G$ be a $p$-group with element breadth 1.  Then $G=(A_1 * A_2 * \cdots * A_k) Z(G)$ where `*' denotes a central product where the isomorphic central subgroups are the derived subgroups of the $A_i$'s.
\end{theorem}
Using these results we will show the following:
\begin{theorem} \label{eb1sb1} Let $G$ be a $p$-group with $\mathrm{ebr}(G)$=$\mathrm{sbr}(G)$=1.  Then $[G:Z(G)]=p^2$ unless $p=2$ and $G$ is isomorphic to one of the following:
\begin{enumerate}
\item $(\Z_2)^n \times Q_8 * D_8,$
\item $(\Z_2)^n \times Q_8 * P_{2,1,1}.$
\end{enumerate}
\end{theorem}
\begin{proof}
Let $G$ be a minimal counterexample to the theorem.  That is, $G$ is a $p$-group with $[G:Z(G)]>p^2$ that is not isomorphic to one of the above groups.  By minimality of $G$ we may assume that $G=A_1 * A_2 \cdots A_k$ where each $A_i$ is minimal non-abelian.  Consider the group $A*B$ where $A$ and $B$ are minimal non-abelian.  Let $z$ be the generator of $G'$ so that $z^p=1$.  We show that we may assume that every non-normal subgroup of either $A$ or $B$ contains $z$.  Suppose not
and suppose that $A$ has a non-normal subgroup $H$ and let $K_1,\cdots,K_{p-1}$ be its distinct conjugates in $A$.  We choose $H$ so that $H$ does not contain $z$.
Let $t$ be a non-central element of $B$ and assume we can pick $t$ such that $z \notin \langle t \rangle$.  Then there is some $b \in B$ such that $[t,b]=z$, therefore $t$ and $tz$
are conjugate.  Clearly the groups $\langle H,t \rangle$, $\langle H,tz \rangle$ and $\langle K_i,t \rangle$ are all distinct and conjugate, a contradiction.  If no such $t$ exists, then every non-normal cyclic subgroup, hence every non-normal subgroup
of $B$ contains $z$.  So either $A$ or $B$ has the given property and we may assume that $p=2$.  Without loss of generality, say $B$ does.
By \ref{blackburnint} we have three possibilities for $B$:

\begin{enumerate}
\item $B \cong Q_8 \times \Z_4 \times (\Z_2)^{n},$
\item $B \cong Q_8 \times Q_8 \times (\Z_2)^n,$
\item $B=\langle g,C \rangle$ with $C$
abelian but not elementary abelian, $1 \ne g^2 \in C$ and
$a^g=a^{-1}$ for all $a \in C$.
\end{enumerate}
The first two are clearly not minimal non-abelian.  By our discussion in the last section, we may assume that $B \cong P_{2,2}$.  Note that $B$, hence $G$, has a normal subgroup $N$ of order 2
besides $\langle z \rangle$ such that $B/N \cong Q_8$.  Therefore consider $G/N$.  By minimality of $G$,
we get that $AN/N \cong D_8$ or $P_{2,1,1}$.  Since $A$ cannot
contain $N$, we get that $A$ must be one of these two groups.  We get the groups [64,201] and
[128,1006] neither of which have subgroup breadth 1 by \ref{64grps} and \ref{128grps}.
\\ \indent Therefore, we may assume that at least one of $A$ or $B$ is Hamiltonian,
so assume that $A$ is. Note that the only minimal non-abelian group is $Q_8$.
A central product of $Q_8$ with $Q_8$ is [32,49] and does not have subgroup breadth 1 by \ref{32grps}, hence
we only must show that $Q_8 * B$ does not have subgroup
breadth 1 when $B \cong P_{i,j}$ or $B \cong
P_{i,1,k}$ when $(i,1,k) \ne (2,1,1)$ or $(1,1,2)$. (We note that
the groups corresponding to (2,1,1) and (1,1,2) are isomorphic.)
Consider the first
case so that $G \cong Q_8 * P_{i,j}$.  We do not use the
standard generators for $Q_8$ as we have already used the variables $i$, $j$ and $k$, therefore we use $x$, $y$ and $z=xy$ in place of them. Suppose first that $i>2$.  Now, $B$ has three
involutions: $a^{2^{i-1}}$, $b^{2^{j-1}}$ and
$a^{2^{i-1}} b^{2^{j-1}}$. Since we require that the common
involution is also a commutator, we must have that $a^{2^{i-1}}=-1$.
If $b$ has order at least 4, then $b^2$ is central, hence we may quotient
out by $\langle b^{2^{j-1}} \rangle$. Hence we may assume that $B$
is a modular group of order at least 16.  So $G$ has a presentation
as follows: $$\langle x,y,a,b \mid
x^4=1,x^2=y^2,[x,y]=x^2=-1,a^{2^i}=1,b^2=1,[a,b]=a^{2^{i-1}},[x,a]=$$
$$[x,b]=[y,a]=[y,b],x^2=a^{2^{i-1}} \rangle.$$  If $i>2$, we claim that the
subgroup $H=\langle b,xa^{2^{i-2}} \rangle$ has more than two
conjugates.  We have $H^y=\langle b,-x a^{2^{i-2}} \rangle$ and
$H^a=\langle ba^{2^{i-1}},xa^{2^{i-2}} \rangle$.  Suppose that
$H=H^y$.  Since both $xa^{2^{i-2}}$ and $-xa^{2^{i-2}}$ are elements
of $H$, we must have that $-1 \in H$. However, $H$ has order 4 and
 $-1$ is not the product of the two generators.  Now suppose
that $H=H^a$.  This says that $a^{2^{i-1}} \in H$, which is also
impossible.  Finally if $H^y=H^a$ then $a^{2^{i-1}} \in H^y$ which
is also impossible.  This proves that we may assume that $i=2$. If $j \ge 2$ we may quotient out by the central
subgroup $\langle b^4 \rangle$ and $G$ has a section isomorphic to
$Q_8 * P_{2,2}$.  This is group 201 of order 64, which
does not have subgroup breadth 1 by \ref{64grps}. This says that
we
may assume that $i=2$ and $j=1$ so that $G=Q_8 * D_8$. \\ \indent
Next we suppose that $B$ is isomorphic to $P_{i,1,k}$. Suppose that both $i$ and $k$ are at
least 2.  Quotienting out by $a^4$ and $c^4$, $G$ must have a
section isomorphic to $Q_8 * P_{2,1,2}$.
This is group 1008 of order 128, so by \ref{128grps} we may assume
that one of $i$ or $k$ is 1.  Note also, that structure-wise, the
roles of $i$ and $k$ are symmetric, hence we might as
well assume that $k=1$. Now suppose that $i>2$.  We quotient out by $\langle a^8 \rangle$
so that $G$ has a section isomorphic to $Q_8 * P_{3,1,1}$; however, this is group 1714 of order 128 which
does not have subgroup breadth 1 by \ref{128grps}.  Hence $G=Q_8 * P_{2,1,1}$.
\\ \indent This shows that only in the two exceptional cases of the statement of the theorem can $k>1$.  We can easily verify that when $k=1$, $[G:Z(G)] \le 4$.  Therefore it remains to show that if $[G:Z(G)]=16$ then $Z(G)$ is elementary abelian.  Suppose that $Z(G)$ has some element $t$ of order 4.  If $t^2 \notin Q_8 * A_2$ where $A_2$ is $D_8$ or $P_{2,1,1}$,
we get groups [128,2162] and [256,26990] respectively, neither of which have subgroup breadth 1 by \ref{128grps} and \ref{256grps}.  Therefore $t^2 \in Q_8 * A_2$.  Specifically, $t^2 \in Z(Q_8 * A_2)$.  When $A_2=D_8$ there is a unique non-identity element of $Z(Q_8 * D_8)$ therefore we only have one choice and we get [64,266] which does not
have subgroup breadth 1 by \ref{64grps}.  When $A_2=P_{2,1,1}$ we get three choices for $t^2$ which result in two different isomorphic classes of groups, [128,2160] and [128,2162], neither of which have subgroup breadth 1 by \ref{128grps}.
\end{proof}
\section{2-Groups with element breadth 2.} \label{ch4}
In this chapter all groups will be 2-groups.
Let $G$ be a group with element breadth 2 and subgroup breadth 1.  We aim to show that $[G:Z(G)] \le 16$.  We first note that by \ref{br2} in a minimal counterexample we must
have that $|G'|=4$.  In this section, $G$ will be a 2-group with element breadth 2 that is a minimal counterexample to \ref{mainresult}.
\begin{lemma} \label{commint} If $G$ has element breadth 2, then $\Omega_1(Z(G)) \le G'$.
\end{lemma}
\begin{proof} Suppose that $z$ is a central involution that is not a
commutator.  Since $G$ is a minimal counterexample, the center
$\overline{Z}$ of $G/\langle z \rangle$ has index at most 16.  We
then have that $[G,Z] \le \langle z \rangle$.  Since $z$ is not a
commutator, we must have that $Z$ is central in $G$, contradicting
that $G$ is a counterexample. \end{proof}
\begin{corollary} $G$ has at most three central involutions.
\end{corollary}
\begin{proof} This follows immediately from $|G'|=4$ \end{proof}
 \indent Now let $z$ be a central involution and let $\pi$ be the natural homomorphism from $G$ to $\overline{G}=G/\langle z \rangle$.
For $g \in G$, $\overline{g} \in Z(\overline{G})$ if and only if $[G,g] \le \langle z \rangle$.  Suppose that there is no $g \in G$ such that $[G,g]=\langle z \rangle$.  Then $\pi^{-1}(Z(\overline{G}))=Z(G)$ and since $z \in Z(G)$ we get that $[G:Z(G)]=[\overline{G}:Z(\overline{G})] \le 16$, a contradiction.  Therefore there is some
element $g \in G$ with $[G,g]=\langle z \rangle$.  This implies the only conjugates of $g$ are $g$ and $gz$, so $C=C_G(g)$ has index 2 in $G$.  Define $$X(g)=\{h \in G \mid [G,h]=\langle t \rangle,C_G(h)=C \}$$
Using the formula that $[a,bc]=[a,c] [a,b]^c$, we see that $H(g)=Z(G) \cup X(g)$ is a subgroup of $C$.  We also see that if $a,b \in X(G)$ then $ab \in Z(G)$.  Therefore $[H(G):Z(G)]=2$.  Since $z \in C$, $H(g)=\pi^{-1}(\overline{C} \cap Z(\overline{G}))$.  This implies that $|Z(G)| \ge \frac{|Z(\overline{G})|}{2}$.  Therefore if $[\overline{G}:Z(\overline{G})] \le 4$ , $[G:Z(G)] \le 16$.  Since $\mathrm{ebr}(\overline{G}$)=1 we get that $\overline{G}$ must be isomorphic to $(\Z_2)^n \times Q_8 * D_8$ or $(\Z_2)^n \times Q_8 * P_{2,1,1}$.  For both of these groups let $A$ be the pre-image under $\pi$ of $(\Z_2)^n$ and let $B$ be the pre-image of the second factor.  We use GAP to determine what $B$ can be.
\begin{theorem} Let $B$ be an extension of the group $Q_8 * D_8$ by $\langle z \rangle$.  Then
$Z(B)$ has index at most 16.
\end{theorem}
\begin{proof}
Using GAP, we obtain the following groups of order 64 that are extensions
of $Q_8 * D_8$:  200, 201,
217, 218, 220, 222, 223, 225, 228, 229, 230, 233, 237, 238, 243,
244, 245, 265.  By \ref{64grps} the only groups from this list that have subgroup breadth 1 are 200, 230, 238, 245, 265 and all of these have center of index at most 16.  \end{proof}
We also note that 200, 238 and 245 are the only ones of these that have no subgroup isomorphic to $D_8$.
\begin{theorem} \label{ext3} Let $B$ be an extension of the group $Q_8 * P_{2,1,1}$ by $\langle z \rangle$.  Then
$Z(B)$ has index at most 16.
\end{theorem}
\begin{proof}
Using GAP, we obtain the following groups of order 128 that are
extensions of $Q_8 * P_{2,1,1}$:  1006, 1008,
1042, 1045, 1048, 1052, 1055, 1059, 1063, 1064, 1068, 1072, 1076,
1083, 1088 1094, 1097, 1103, 1110, 1113, 1114, 1714, 1715, 1716,
1717, 1718, 1719, 2158.  By \ref{128grps} the only groups from this list that have subgroup breadth 1 are 1716 and 2158 and both of these have center of index at most 16.  \end{proof}
Now, since $\overline{A}'=\overline{1}$ we have that $A$ must have element breadth at most 1.  Next we consider the possible structures of $A$.
\\ \indent By \ref{eb1sb1}, if $A$ is non-abelian, $A$ can either be written as $(\Z_2)^n \times (Q_8 * D_8),(\Z_2)^n \times (Q_8 * P_{2,1,1})$ or as $CZ(A)$ where $C$ is minimal non-abelian.  However $\overline{A}$ must be elementary abelian.  The first case contains dihedral subgroups, which is impossible by \ref{johnd8}.  In the second case $\overline{A}$ is clearly not elementary abelian as $A$ has rank 4.  For the third case, $Q_8$ is the only possibility for $C$ such that $\overline{A}$ is elementary abelian.  Now $Z(A)$ can have elements of order at most 4. Let $t \in Z(A)$ such that $t^4=1$.  If $t^4 \in Q_8$ it is easily verified that $\langle Q_8,t \rangle \cong Q_8 \times \Z_2$.  As $z \in Q_8$ we may, therefore assume that $A \cong Q_8 \times D$ where $D$ is elementary abelian or that $A$ is abelian with $\overline{A}$ elementary abelian.  Consider the first case.  Let $t \in D$.  Then the group $tB$ has a center of index at most 16.  Also, by \ref{commint}, $t$ is not central.  We note that $Z(tB)=C_{B}(t) \cap Z(B)$.  For all of the possible groups $B$ except [128,2158], $Z(B) \le \Phi(B)$.  As $C_B(t)$ is maximal in $B$, we have $Z(tB)=Z(B)$.  Therefore $tB$ has a center of index 32, a contradiction.  We may, therefore assume that $A \cong Q_8$ in this case.  If $A$ is abelian, the same argument shows that $A \cong \Z_4$ or $\Z_2$.  This shows that $|G| \le 512$.  By \ref{32grps}, \ref{64grps}, \ref{128grps} and \ref{256grps}, we may assume that $|G|=512$ so that $A \cong Q_8$ and $B$ is either [128,1716] or [128,2158].
\\ \indent Both contain a subgroup isomorphic to $Q_8$, therefore if $[A,B]=1$ then $G$ contains a subgroup isomorphic to either $Q_8 * Q_8$ or $Q_8 \times Q_8$, neither of which have subgroup breadth 1.  Therefore $[A,B]=\langle z \rangle$.  If $B$=[128,1716] we can verify in GAP that this group has the structure
$(\Z_8 \times Q_8) \rtimes \Z_2$ and presentation:
$$\langle a,b,c,d \mid a^8=b^4=d^2=1,c^2=[b,c]=b^2,$$ $$[a,b]=[a,c]=[b,d]=[c,d]=1,[a,d]=a^4 b^2 \rangle .$$
In this presentation $a^4=z$.  Therefore, we can find some group of order 256 of the form $A(\Z_8 \times Q_8)$.  Using \ref{256grps} we can check that only groups 6648, 26461, 26462, 53175 and 53232 have element breadth 2, subgroup breadth 1 and normal subgroup isomorphic to $\Z_8 \times Q_8$ and none of these groups has a normal subgroup isomorphic to $Q_8$. (Note that a group with element breadth 1 cannot have this structure by \ref{eb1sb1}.)
\\ \indent Therefore, we may assume that $B$ is [128,2158] which is isomorphic to $\Z_2 \times Q_8 * ((\Z_4 \times \Z_2) \rtimes \Z_2)$.  Consider the subgroup $iB$.  It is easily shown that if $C=(\Z_4 \times \Z_2) \rtimes \Z_2$ that $iB \cong \Z_4C$.  We verify in GAP that no group in \ref{256grps} has this structure.  This completes the proof for 2-groups.
\section{$p$-groups with element breadth 2 for odd $p$}
Our proof when $p$ is an odd prime will differ greatly from that for
$p=2$, both in format and in length.  There are several underlying
reasons for this disparity.  It can be shown that every
irreducible representation of a 2-group with subgroup breadth 1 has
degree at most 4.  The tensor product of the respective  2-dimensional irreducible representations of $Q_8$ and
$D_8$ shows that this bound is sharp for $Q_8 * D_8$.  Much can
be said about groups all of whose irreducible representations have degree dividing
$p^2$; however a proof using these facts would likely be fairly complex.
By a theorem of Isaacs, much more can be said when all non-linear irreducible representations
have degree $p$ (regardless of parity):
\begin{theorem} \cite{rx1} \label{rx1} (Theorem 12.11) A $p$-group, $G$, has irreducible representations of only
degrees 1 and $p$ if and only if one of the following holds:
\begin{enumerate}
\item $G$ has a maximal subgroup which is abelian,
\item $[G:Z(G)]=p^3$.
\end{enumerate}
\end{theorem}
By adding a hypothesis, we can show that (1) implies (2).
\begin{theorem} \label{ref2} Let $G$ be a $p$-group with subgroup breadth 1 and element
breadth 2 that is a minimal counterexample to $[G:Z(G)] \le p^3$.  Then $G$
does not have an abelian maximal subgroup.
\end{theorem}
\begin{proof} Let $|G|=p^n$.
By Lemma 12.12 in \cite{isaacs}, if $A$ is an abelian subgroup with
$G/A$ cyclic, we have $|A|=|G'||A \cap Z(G)|$.  By \ref{br2} we may assume
that $|G'|=p^2$.  Then we have $[A:A \cap Z(G)]=p^2$.  This clearly implies that
$Z(G)$ has index at most $p^3$.  Therefore $G$ cannot have an abelian maximal subgroup.
\end{proof}
We now use a series of smaller results to prove our main theorem:
 \begin{theorem} \label{oddsb1} Let $G$ be a $p$-group with subgroup breadth 1, where $p$ is an odd prime.  Then
$[G:Z(G)] \le p^3$.
\end{theorem}
\begin{proof} Let $G$ be a minimal counterexample to the claim that $\mathrm{sbr}(G)=1$ and $[G:Z(G)] \le p^3$.  By \ref{rx1} and \ref{ref2}, $G$ is also a minimal counterexample to the claim that every non-linear irreducible representation of a group with subgroup breadth 1 has degree at most $p$.  Let $\phi$ be a representation of degree $p^i$ where $i>1$. By $\ref{eb1sb1}$, we may assume that $G$ has element breadth 2, and by \ref{br2} we may assume that $|G'|=p^2$.
\begin{enumerate}
\item $Z(G)$ is cyclic.
 \begin{proof} Suppose not.  By Schur's lemma, the image of any irreducible representation has a cyclic center, therefore $\phi$ is not faithful.  Hence, $G/\text{ker}(\phi)$ also has an irreducible representation of degree $p^i$; however all quotient groups of $G$ must also have subgroup breadth 1, therefore this is a smaller counterexample.
\end{proof}
\item Let $X=\langle x \rangle$ be a subgroup of $G$ of order $p$.  Then $X$ commutes with all of its $G$-conjugates.
\begin{proof} Since $N_G(X)/C_G(X)$ is a subgroup of $\text{Aut}(X) \cong \Z_{p-1}$, we have $N_G(X)=C_G(X)$.  Also for $g \in G$, we have $$N_G(X^g)=N_G(X)^g=N_G(X)=C_G(X)=C_G(x)$$ and clearly $x^g$ normalizes $X^g$. (The second equality follows from the fact that $N_G(X)$ is a maximal subgroup of a nilpotent group, hence is normal.)  Therefore, either $X \lhd G$ or $N_G(X)$ is a maximal subgroup of $G$.
\end{proof}
\item If $Y=\langle y \rangle$ is another subgroup of $G$ of order $p$ then $[X,Y]=1$.
\begin{proof}
Suppose not.  Then $X$ has $p$ $Y$-conjugates, call them $X_1, \cdots, X_p$.  Since $G$ has subgroup breadth 1, these are all of the $G$-conjugates of $X$.  Hence $N=\langle X_1,\cdots,X_p \rangle \lhd G$.  By (2), $N$ is elementary abelian, of order $p^k$ for some $k>0$.  Let $C=C_G(N)$.  Since $$C_G(X^g)=N_G(X^g)=N_G(X)^g=N_G(X)=C_G(X)$$ we have that $C=C_G(X)$.  Since $[X,Y] \ne 1$, we have $Y$ is not contained in $C$, hence $G=CY$.  Therefore, any element of $N$ that centralizes $Y$ is central in $G$.  However, since $Z(G)$ is cyclic, we must have that $C_N(Y)$ has order $p$.  Since $NY$ is a semidirect product with kernel $N$, we have that $N_N(Y)=C_N(Y)$.  Since $G$ has subgroup breadth 1, we have $[N:N_N(Y)] \le p$, hence we have $k \le 2$.  If $k=1$, then $[X,Y]=1$, a contradiction.  Hence assume $k=2$.  Then $NY=\langle x,y \rangle$ is an extraspecial group of order $p^3$. Now, both $X$ and $Y$ have $p$ conjugates in $NY$, hence $NY \lhd G$.  We claim that $$[G,YN]=[YN,YN]=Z(YN).$$  We have $G=YC$.  Then let $y_1,y_2 \in Y$, $c \in C$ and $n \in N$. Then
$$[y_1 c,y_2 n]=[y_1 c,n][y_1 c,y_2]^n=[y_1,n]^c [c,n] [y_1,y_2]^{nc} [c,y_2]^n=[c,n] [c,y_2]^n$$
Now, since $YN \lhd G$, we have $[YN,YN] \lhd G$ hence we only must show that if $c \in C$ and $y \in Y$, then $[c,y_2] \in [YN,YN]$.  Now, since $C_G(y_2)$ has index $p$ in $G$ and $C_C(y_2)$ has index $p$ in $C$ we have that $D=C_C(y_2)$ has index $p$.  Then $C=DX$.  We therefore have $$[c,y]=[dx_1,y_2]=[d,y_2]^{x_1} [x_1,y_2]=[x_1,y_2]$$ which is clearly a commutator of $YN$.  So $[G,YN] \le [YN,YN]$.  Clearly we have the other inclusion, so $[G,YN]=[YN,YN]$.  Since $YN$ is extraspecial, we have that $[G,YN]=Z(YN)$.  Let $H=C_G(YN)$.  Then $[H,YN]=1$ and $H \cap YN=Z(YN)$.  Also, since $H=C_G(Y) \cap C_G(N)$ has index at most $p^2$ and $|YN|=p^3$ we get $G=YNH$.  Therefore, $G$ is a central product of $YN$ with $H$.  Now, since $YN$ has non-normal subgroups if $H$ has non-normal subgroups, they must all contain $Z(YN)$. However, by \ref{blackburnint} this only is possible when $p=2$.  Therefore, $H$ is abelian.  Since $X \nleq H$, and $X \le YN$, we have $HX$ is also abelian, and, having index $p$, produces a contradiction.
\end{proof}
\item $\Omega_1(G)$ is an elementary abelian group of size $p^3$.
\begin{proof}
 Let $O=\Omega_1(G)$.  By (3), $O$ is elementary abelian.  Let $Z=Z(G) \cap O$.  Since $Z(G)$ is cyclic, $|Z|=p$.  If $O=Z$ then $G$ has a unique subgroup of order $p$, hence is cyclic.  So we may assume that $Z \ne O$.  Let $X$ and $Y$ be subgroups of $O$ of order $p$ both not equal to $Z$.  If $N_G(X) \ne N_G(Y)$ then $X \times Y$ has at least $p^2$ conjugates (neither $X$ nor $Y$ can be normal because then $X$ or $Y$ is central).  Therefore, all non-central subgroups of $G$ have the same centralizer, so $C=C_G(O)=C_G(X)$, so $[G:C]=p$.  Therefore, $G$ acts on $O$ as a linear group $L$ of order $p$.  Let $g$ generate $L$.  If $g$ is in Jordan form with respect to some basis for $O$ then, since $C_O(g)$ is contained in $Z(G)$, we have $C_O(G)=Z$. Therefore, $g$ has exactly one Jordan block.  Also, this means that if $x \in G-C$ then $C_O(x)=Z$.  Since $G$ has element breadth 2, we have $[O:Z] \le p^2$, hence $|O| \le p^3$. By a result of \cite{nearuniform}, if $|\Omega_1(G)|=p^2$ then $G$ is either metacyclic or a 3-group of maximal class.  If $G$ is a metacyclic group, \cite{csw} says that $G$ must have element breadth 1; however, our classification of such groups shows that this implies $[G:Z(G)]=p^2$.  If $G$ is a 3-group of maximal class, this implies that $[G:G']=p^2$ which implies that $G$ has order $p^4$.  Since $p$-groups have non-trivial centers, we clearly have $[G:Z(G)] \le p^3$, a contradiction.  Therefore, $O$ has order $p^3$.
\end{proof}
\end{enumerate}
We now complete the proof of the theorem.  Let $C=C_G(O)$ and let $h \in G-C$.  Let $H=\langle h \rangle$ and $N=N_G(H)$.  As above, $h$ acts with one Jordan block on $O$, therefore we can find $x \in O$ such that $x \notin N$.  Let $X=\langle x \rangle$.  Therefore, $G$ is a semidirect product of $X$ and $N$, and we have $|\Omega_1(N)|=p^2$. (If $|\Omega_1(N)|=p$ we get that $N$ is cyclic, therefore $G$ is metacyclic and cannot have element breadth 2.)  By the result above, $N$ is either metacyclic or a 3-group of maximal class.  If $N$ is metacyclic by \cite{csw} $N$ has element breadth 1 and by \ref{eb1sb1} we have that $N=\langle a,b \mid a^{p^{n-1}}=b^p=1,[a,b]=a^{p^{n-2}} \rangle$. Therefore, $\Phi(N)=Z(N)$.  Now, since $\mathrm{ebr}(N)$=1, we have that $C_N(H)$ has index $p$, hence is maximal and must therefore contain $Z(N)$.  This shows that $Z(N)=Z(G)$.  Hence $[G:Z(G)]=[G:Z(N)]=[G:N][N:Z(N)]=p^3$.  If $N$ is a 3-group of maximal class, then, as above, we have that $N$ has order $3^4$, hence $G$ has order $3^5$.  Now, $N_G(O)=C_G(O)$ which has index $3$.  Clearly $O \le Z(C_G(O))$ and since a $p$-group cannot have a center of index $p$, we get that $C_G(O)$ must be an abelian subgroup. (Note that this argument does not produce an abelian maximal subgroup in the generic case since we have no guarantee that all elements of $C_G(O)$ must commute with each other.)  This final contradiction completes the proof.
\end{proof}
We note that the proof of \ref{oddsb1} does not necessarily imply that every $p$-group with element breadth 2 and subgroup breadth 1 has an abelian maximal subgroup since we were only examining the structure of a minimal counterexample to the theorem. However, we have no examples of $p$-groups for odd primes $p$ with subgroup breadth 1 which do not have an abelian maximal subgroup.
\section{$p$-groups with subgroup breadth more than one}
We conclude with a conjecture regarding $p$-groups with subgroup breadth more than one:
\begin{conjecture}
Let $G$ be a $p$-group with subgroup breadth $k$. Then \begin{equation*}
|G:Z(G)| \le \begin{cases} 2^{3k+1} \text{ if } p=2, \\ p^{3k} \text{ if } p>2. \end{cases} \end{equation*}
\end{conjecture}
It can be verified in GAP that there exist no groups of order $2^9$ with subgroup breadth 2 and center of order 2 and similarly there exist no groups of order $3^8$ with subgroup breadth 2 and center of order $3$.  It should be mentioned that both of these are substantial computations.  The groups of order 512 are in the Small Groups Library, however, there are 10,494,213 such groups.  Using trivial parallelization, approximately 50 processors were used to get a list of the sizes of the centers of all groups of order 512.  There are 5,327 groups with center of size 2, only 10 of which have cyclic breadth 2.  We rule these out case-by-case.  The groups of order $3^8$ were obtained using the $p$-group generation algorithm \cite{pquotient} implemented in the GAP package ANUPQ.
\\ \indent It is clear that the methods of the proof of \ref{mainresult} in this paper cannot be extended to even the case $k=2$.  It should be noted that classifications of $p$-groups with element breadth 3 exist (see, for instance, \cite{ps} or \cite{2br3}).  However, no generic classification of $p$-groups with element breadth 4 exist.
\section{Computational Results}
In this section we provide a complete list of all 2-groups which have subgroup
breadth 1, up to order 256.  These calculations were all made in GAP.
\begin{theorem} \label{32grps} A non-abelian group of order 32 has subgroup breadth at most 1 if and
only if it is one of the following:
\\ 2, 4, 5, 8, 12, 15, 17, 22, 23, 24, 25, 26, 29, 32, 35, 37, 38, 41, 46, 47, 48,
50. \end{theorem}
\begin{theorem} \label{64grps} A non-abelian group of order 64 has subgroup breadth at most 1 if and
only if it is one of the following:
\\ 3, 17, 22, 27, 29, 44, 45, 51, 56, 57, 58, 59, 84, 85, 86, 87, 88, 93, 103, 104, 105,
110, 112, 113, 115, 126, 127, 184, 185, 193, 194, 195, 196, 197,
198, 200, 204, 208, 212, 214, 230, 238, 245, 247, 248, 252, 261,
262, 263, 265. \end{theorem}
\begin{theorem} \label{128grps} A non-abelian group of order 128 has subgroup breadth at most 1 if and
only if it is one of the following: \\ 5, 43, 44, 106, 108, 129, 131,
153, 154, 160, 164, 180, 181, 182, 183, 184,
  457, 458, 459, 460, 469, 476, 477, 480, 481, 483, 498, 499, 501, 509, 838,
  839, 840, 843, 844, 881, 882, 883, 884, 894, 895, 899, 914, 915, 989, 990,
  998, 999, 1000, 1001, 1002, 1003, 1004, 1602, 1603, 1604, 1606, 1608, 1609,
  1618, 1634, 1635, 1636, 1646, 1649, 1650, 1652, 1658, 1690, 1691, 1692,
  1696, 1716, 2137, 2138, 2151, 2152, 2153, 2154, 2155, 2156, 2158, 2165,
  2169, 2173, 2175, 2198, 2208, 2262, 2302, 2303, 2308, 2320, 2321, 2322,
  2324. \end{theorem}
\begin{theorem} \label{256grps} A non-abelian group of order 256 has subgroup breadth at most 1 if and
only if it is one of the following: \\  40, 124, 126, 317, 319, 453,
455, 498, 500, 531, 532, 538, 827, 828, 829,
  830, 831, 835, 1119, 1131, 1247, 3680, 3681, 3683, 3692, 4385, 4386, 4387, 4388, 4389,
  4390, 4395, 4396, 4398, 4399, 5526, 5527, 5531, 5532, 5536, 5578, 5579, 5586, 5587, 5598, 5640, 5641,
  5643, 5649, 6535, 6536, 6537, 6540, 6541, 6614, 6615, 6616, 6617, 6627,
  6628, 6632, 6647, 6648, 6724, 6725, 6733, 6734, 6735, 6736, 10299, 10300, 10301, 10303, 10305, 10306, 10311, 13314, 13315, 13316,
  13318, 13320, 13321, 13337, 13349, 13350, 13352, 13359, 13361, 13364,
  13368, 13411, 13413, 13416, 13420, 13444, 13449, 13451, 13454, 13459,
  13471, 13503, 26309, 26310, 26311, 26313, 26318, 26319, 26382, 26383, 26384, 26387,
  26403, 26404, 26406, 26417, 26460, 26461, 26462, 26466, 26960, 26961,
  26974, 26975, 26976, 26977, 26978, 26979, 26980, 26981, 53039, 53040, 53041, 53043, 53047, 53048, 53060, 53082, 53083, 53084,
  53098, 53102, 53103, 53105, 53118, 53173, 53174, 53175, 53184, 53232,
  55609, 55610, 55627, 55628, 55629, 55630, 55631, 55632, 55634, 55645, 55649, 55653, 55655, 55680, 55684, 55697, 56060, 56061, 56067, 56083, 56084, 56085, 56087.
  \end{theorem}
\bibliographystyle{plain}	
\bibliography{myrefs}		
\end{document}